\documentclass[twoside,leqno]{article}
\usepackage{amsmath}
\newcommand\ack{\section*{Acknowledgement}}

\newcommand\ds{\displaystyle}
\newcommand{\Real}{\hbox{\bf R}}

\newtheorem{remark}{\bf Remark}[section]

\newtheorem{theorem}{\bf Theorem}[section]

\newenvironment{proof}{
\begin{trivlist}
\item[\hspace{\labelsep}{\it\noindent Proof: }]
}{\par\hfill\fine\end{trivlist}
\par}
\input{epsf.sty}
\newenvironment{figMacPc}[4]				
{\begin{figure}[htb]						
	\epsfxsize=#2\centerline{\epsfbox{#1}}		
	\caption{#3}					
	\label{#4}							
}
{\end{figure}}
\newcommand{\fine}{\quad\rule{2mm}{2mm}}		    % box pieno
\begin{document}
\begin{center}
{\Large \bf Towards Dead Time Inclusion  \\ \medskip in Neuronal Modeling}
\end{center}
\par\noindent
\begin{center}
{\bf \large A. Buonocore$^{(1)}$, G. Esposito$^{(1)}$, V. Giorno$^{(2)}$\\
 and C. Valerio$^{(1)}$} 
\end{center}  
\par\noindent
\begin{center}
{\it (1)} Dipartimento di Matematica e Applicazioni, Universit\`a di Napoli Federico II, 
Via Cintia, Napoli, Italy, Email: \{aniello.buonocore\}@unina.it\\
{\it (2)} Dipartimento di Matematica e Informatica, Universit\`a di Salerno, Via Allende,  
Baronissi (SA), Italy, Email: \{giorno\}@unisa.it
\end{center}
\par\noindent
\begin{abstract} 
A mathematical description of the refractoriness period 
in neuronal diffusion modeling is given and its  moments  
are explicitly  obtained in a form that is suitable for quantitative 
evaluations. Then, for the Wiener, Ornstein-Uhlenbeck and 
Feller neuronal models, an analysis of the features exhibited by the mean and variance of 
the first passage time and of refractoriness period is performed.
\end{abstract} 
%
%
%\key{Refractoriness; Neuronal models; Output distribution; Diffusion models} 

%\subj{60J60; 60J70; 92C20} 
%------------------------------------------------------
\section{Introduction}
%------------------------------------------------------
\hfill\\
Mathematical descriptions of Òdead timeÓ, or refractoriness, in neuronal modeling 
have long traditions dating back at least to mid sixties when special attention was devoted to the 
description of the evolution of networks of switching elements whose behavior was meant 
to simulate that of physiological neurons via  certain rather drastic 
simplifications \cite{Ricciardi66}. Furthermore, the approach to neural modeling was shown 
to bear certain strong analogies with the stochastic description of the input-output 
features of radioactive particle counters. In such context, as early as 1948, W. Feller 
proved that under a suitable formulation all problems concerning single counters 
reduce to special instances of the theory of the summation of random variables. 
Exploiting the above mentioned analogies, the simplest neural model may be conceived as a 
black box possessing the following distinctive features: {\it (i)} it is a 
threshold element, {\it (ii)} its output response consists of pulses of constant 
amplitude and width and {\it (iii)} there exists a constant dead time. More 
accurately, one could define this dead time also as a deterministic 
function of certain measurable parameters, such as time or input pulse 
amplitude, or view it as a stochastic process. 
\par
As a first attempt towards 
a quantitative treatment of the dead time effects in neural modeling, we 
look at the input of the neuron as a randomly distributed Poisson-type pulse 
train. Its output is then determined by imposing the restriction that 
following each input pulse a dead time period is activated during which no 
further pulses can be produced at the output. Even for such an oversimplified 
instance the investigation of the role played by the 
dead time in determining the distribution of the output when the input is described 
by a given distribution is a very challenging task. 
\par
Let $\tau$ denote this dead time, 
i.e. the time interval following every firing 
during which the neuron cannot fire again. Let us assume that the net 
input to the neuron in time interval $(0,T)$ is modeled 
by a sequence of positive pulses of standard strength whose time of occurrences are 
Poisson distributed with rate $\lambda>0$. We purpose to determine the distribution 
$\Pi_n(T,\tau)$ of the output pulses as a function of dead time $\tau$. A rather 
cumbersome amount of computations leads one to conclude that the assumed input 
distribution
\begin{equation}
P_n(T)={(\lambda\,T)^n\over n!}\;e^{-\lambda\,T}, \quad T>0,\;n=0,1,2,\ldots
\label{eq:Poisson}
\end{equation}
generates the following firing distribution valid for all $n\geq 1$ (cf. \cite{Ricciardi66}) :
\begin{eqnarray}
&&\Pi_n(T,\tau)=\vartheta[T-(n-1)\tau]\;\biggl\{1-e^{-\lambda\,[T-(n-1)\tau]}
\sum_{k=0}^{n-1}{\lambda^k\,[T-(n-1)\tau]^k\over k!}\biggl\}\nonumber\\
&&\hspace*{1.6cm} -\vartheta(T-n\tau) \biggl[1-e^{-\lambda\,(T-n\tau)}
\sum_{k=0}^n{\lambda^k\,(T-n\tau)^k\over k!}\biggr],
\label{eq:FiringDistrib}
\end{eqnarray}
where $\vartheta(x)$ denotes the  Heaviside unit step function:
\begin{equation}
\vartheta(x)=\left\{
\begin{array}{ll}
1,&x>0\\
0,&x\leq 0.
\end{array} \right.
\label{theta}
\end{equation}
Although the stated problem 
has been the object of several investigations (see, for instance, \cite{Teich78}), a quantitative 
evaluation of the effect of dead time on the statistical parameters of the 
output appears to be still lacking. 
\par
In the remaining part of this paper, we shall outline a totally different approach towards 
the inclusion of refractoriness in the neuronal model. As in  \cite{Buonocore02one} 
and \cite{Buonocore02two}, we model the time course of the membrane potential 
by a time-homogeneous diffusion process and then 
assume that the firing threshold acts as some kind of elastic boundary characterized 
by preassigned reflection and absorption parameters. In other words, 
we assume that an action potential is released whenever the process first attains 
the firing threshold. After the firing, a period of refractoriness of random 
duration occurs, at the end of which the process is instantaneously reset 
at a fixed state. Then, the subsequent evolution of the action potential proceeds as before, 
until the threshold is again reached. A new firing then occurs, followed by a 
new period of refractoriness, and so on. Use of the above approach allows one to 
mimic the effects of refractoriness for the specified neuronal model. 
\par 
In order to be able to apply the specified paradigm to the description of neuronal 
models in the presence of refractoriness, an investigation of certain general 
features of diffusion processes in the presence of an elastic boundary is 
necessary. This task will be accomplished in  
Section 2, where we shall analyze the features of the moments of the random 
variable modeling the neuron's intrinsic refractoriness. In Section 3, a specific analysis 
will be provided of three neuronal models based on the Wiener, Ornstein-Uhlenbeck and Feller diffusion 
processes, and a comparative discussion of the refractoriness features exhibited by these models will be 
performed.
%
%------------------------------------------------------
\section{Effect of Refractoriness}
%------------------------------------------------------
\hfill\\
Let $\{X(t),t\geq 0\}$ be a regular, time-homogeneous diffusion process,  
defined over the interval $I=(r_1,r_2)$, characterized by 
drift and infinitesimal variance $A_1(x)$ and $A_2(x)$, respectively. 
Throughout, we shall assume that Feller conditions on these functions are fulfilled 
\cite{Feller52}. 
Let $h(x)$ and $k(x)$ denote scale function and speed density of $X(t)$:
$$
h(x)=\exp\biggl\{-2\int^x{A_1(z)\over A_2(z)}\;dz\biggr\},\qquad
k(x)={2\over A_2(x)\,h(x)}
$$
and  
$$
H(r_1,y]=\int_{r_1}^y  h(z)\;dz,\qquad K(r_1,y]=\int_{r_1}^y  k(z)\;dz$$
scale and  speed measures, respectively. 
\par
We define the random variable \lq\lq first passage time\rq\rq\ (FPT) of $X(t)$ 
through $S$ $(S\in I)$ with $X(0)=x<S$:
\begin{equation}
T_x=\inf_{t\geq 0}\{t:X(t)\geq S\},\qquad X(0)=x<S.
\label{eq:FPTrandomvariable}
\end{equation}
Then,
\begin{equation}
g(S,t|x)={\partial\over\partial t} P(T<t),\qquad x<S
\label{eq:FPTdensity}
\end{equation}
is the FPT pdf of $X(t)$ through $S$ conditional upon $X(0)=x$. 
\par
In the neuronal modeling context the state $S$ represents the neuron's firing 
threshold and $g(S,t|x)$ the firing pdf. 
\par
More realistically then in past approaches, here we shall assume that after 
each firing a period of refractoriness  of random duration occurs, during 
which either the neuron is completely unable to respond, or it only partially responds  
to the received stimulations. To this end, we look at the threshold $S$  
as an elastic barrier being \lq partially transparent\rq, in the sense that 
its behavior is intermediate between total absorption and total reflection. 
The degree of elasticity of the boundary depends on the choice of two parameters, 
$\alpha$ (absorbing coefficient) and $\beta$ (reflecting coefficient), 
with $\alpha>0$ and $\beta\geq 0$. Hence, $p_{_{R}}:=\beta/(\alpha+\beta)$ denotes the 
reflecting probability at the boundary $S$, and $1-p_{_{R}}=\alpha/(\alpha+\beta)$ the 
absorption probability at $S$. 
%
%FIGURE1 REFRACTORINESS
\begin{figMacPc}
			{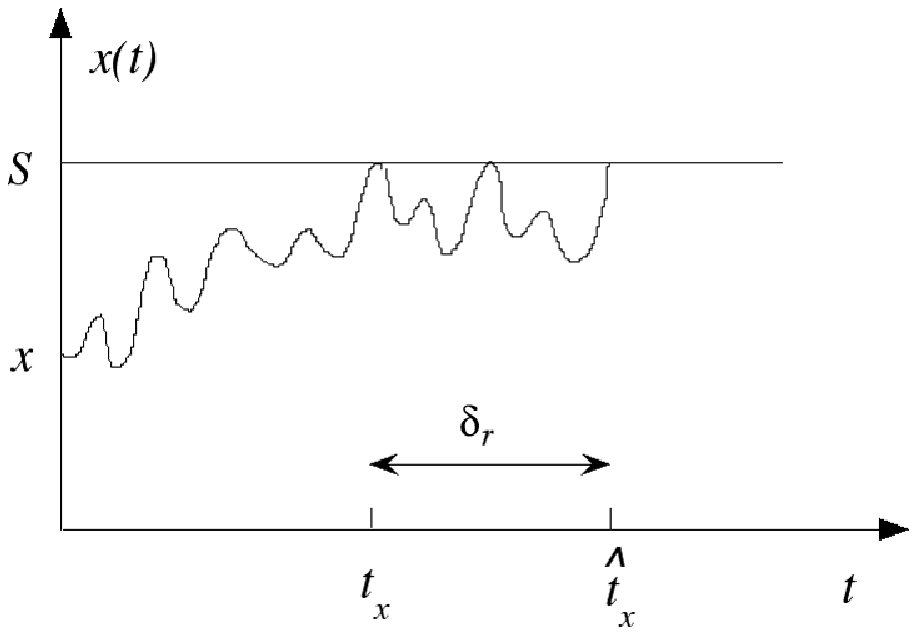}
			{7 cm}
                  {\small Illustrating first passage time $t_x$ through $S$ (i.e. firing 			time), neuron's refractoriness period  $\delta_r$ and first exit time  			$\widehat{t}_x$ for a single sample path $x(t)$ of $X(t)$. By $t_x$, 			$\widehat{t}_x$ and $\delta_r$ we have indicated the appropriate values of 			$T_x$, $\widehat{T}_x$ and $T_r$, respectively.}
			{fig:ref1}
\end{figMacPc}
%\begin{figure}[t]
%\includegraphics{fig1.ps}
%\caption{\small Illustrating first passage time $t_x$ through $S$ (i.e. firing time), 
%neuron's refractoriness period  $\delta_r$ and first exit time  $\widehat{t}_x$ for a single 
%sample path $x(t)$ of $X(t)$. By $t_x$, $\widehat{t}_x$ and $\delta_r$ we have indicated 
%the appropriate values of $T_x$, $\widehat{T}_x$ and $T_r$, respectively.}
%\label{fig:ref1}
%\end{figure}
%
We denote by $\widehat{T}_x$ the random variable describing the 
\lq\lq first exit time\rq\rq\ (FET)  of $X(t)$ through $S$  if 
$X(0)=x< S$, and by $g_e(S,t|x)$ its pdf. The random variable 
$T_r$ will denote the \lq\lq refractoriness period\rq\rq\ and $g_r(S,t|S)$ 
its pdf. Since $\widehat{T}_x$ can be viewed as the sum of random 
variable $T_x$ describing the first passage time through $S$ (firing time) 
and of $T_r$ (see Figure \ref{fig:ref1}) one has:
\begin{equation}
g_e(S,t|x)=\int_0^tg(S,\tau|x)\,g_r(S,t|S,\tau)\;d\tau.
\label{eq:integralequat}
\end{equation}
In the sequel we assume that one of the following cases holds:
\begin{description}
\item{\em (i)} $r_1$ is a natural nonattracting boundary  and $K(r_1,y]<+\infty$;
\item{\em (ii)} $r_1$ is a reflecting boundary or it is an entrance boundary. 
\end{description}
Under such assumptions, if $x<S$ the first passage probability $P(S|x)$ from 
$x$ to $S$ is unity and the FPT moments $t_n(S|x_0)\equiv E(T_x^n)$ are finite and can be 
iteratively calculated as 
\begin{eqnarray}
&&t_n(S|x):=\int_0^{\infty}t^n\,g(S,t|x)\;dt=n\,\int_x^Sh(z)\;dz \int_{r_1}^zk(u)\,t_{n-1}(S|u)\;du
\quad (n=1,2,\ldots),
\label{eq:momentsFPT}
\end{eqnarray}
where $t_0(S|x)=P(S|x)=1$ (cf., for instance, \cite{Siegert51}).
%***************
% THEOREM 2.1
%***************
\begin{theorem}
Under the assumption {\it (i)} and {\it (ii)},  if $\alpha>0$ the first exit time probability 
$$
{\widehat P}(S|x):=\int_0^{+\infty}g_e(S,t|x)\;dt\qquad  (x< S)
$$ 
is unity.
\end{theorem}
\begin{proof}
We consider separately the cases {\it (i)} and {\it (ii)}.
\par\noindent
{\it Case (i)}\quad Let $\widehat{P}(S_1,S|x)$ $(r_1<S_1<x<S)$ be the first exit time probability 
through the elastic boundary $S$ in the presence of an absorbing boundary $S_1$. This  
is solution of the differential equation
\begin{equation}
A_1(x)\,{d\psi_0(x)\over dx}+{A_2(x)\over 2}\,{d^2\psi_0(x)\over dx^2}=0
\label{eq:diffequationprob}
\end{equation}
subject to conditions
\begin{equation}
\lim_{x\downarrow S_1} \psi_0(x)=0,\qquad
 \alpha\,\lim_{x\uparrow S} \bigl[1-\psi_0(x)\bigr]-\beta\,\lim_{x\uparrow S}
\Bigl\{h^{-1}(x)\,{d\psi_0(x)\over dx}\Bigr\}=0.
\label{eq:diffequationprobcond1}
\end{equation}
Since
\begin{equation}
A_1(x)\,{d\psi_0(x)\over dx}+{A_2(x)\over 2}\,{d^2\psi_0(x)\over dx^2}\equiv
{1\over k(x)}\;{d\over dx}\biggl[{1\over h(x)}\,{d\psi_0(x)\over dx}\biggr]
\label{eq:diffequationprobrelation}
\end{equation}
from (\ref{eq:diffequationprob}) one has
\begin{equation}
\psi_0(x)=A+B\int^xh(z)\;dz,
\label{eq:gensolution}
\end{equation}
where $A$ and $B$ are arbitrary real constants. By imposing boundary conditions 
(\ref{eq:diffequationprobcond1}), one obtains
\begin{equation}
\widehat{P}(S_1,S|x)={\alpha\,{\ds\int_{S_1}^xh(z)\;dz}\over \alpha\,{\ds\int_{S_1}^Sh(z)\;dz}+\beta}\cdot
\label{eq:diffequationprobgensolut}
\end{equation}
Since $r_1$ is a natural nonattracting boundary one has $H(r_1,x]=+\infty$; hence, making use of  
(\ref{eq:diffequationprobgensolut}), one has
$$
\widehat{P}(S|x):=\lim_{S_1\downarrow r_1}\widehat{P}(S_1,S|x)=1,
$$
where the last equality follows by exploiting l'Hospital's rule.  
\par\noindent
{\it Case (ii)}\quad Let $\widehat{P}(S|x)$ $(r_1<x<S)$ the first exit time probability 
through the elastic boundary  $S$ in the presence of a reflecting boundary 
or of an entrance boundary $r_1$. This is solution of the differential equation 
(\ref{eq:diffequationprob}) subject to conditions:
\begin{equation}
\lim_{x\downarrow r_1} \Bigl\{h^{-1}(x)\,{d\psi_0(x)\over dx}\Bigr\}=0,\qquad
\alpha\,\lim_{x\uparrow S} \bigl[1-\psi_0(x)\bigr]-\beta\,\lim_{x\uparrow S}
\Bigl\{h^{-1}(x)\,{d\psi_0(x)\over dx}\Bigr\}=0.
\label{eq:diffequationprobcond2}
\end{equation}
Since (\ref{eq:diffequationprobrelation}) holds, from (\ref{eq:diffequationprob}) 
one obtains again the general solution (\ref{eq:gensolution}). By imposing  boundary conditions 
(\ref{eq:diffequationprobcond2}), one finally has $\widehat{P}(S|x)=1$.
\end{proof}
%
%***************
% THEOREM 2.2
%***************
\begin{theorem}
Under the assumption {\it (i)} and {\it (ii)},  if $\alpha>0$   the first exit time moments 
${\widehat t}_n(S|x)\equiv E({\widehat T}_x^n)$ can be iteratively calculated as
\begin{eqnarray}
&&\widehat{t}_n(S|x):=\int_0^{\infty}t^n\,g_e(S,t|x)\;dt=n\,\Biggl\{\int_x^Sh(z)\,dz
\int_{r_1}^zk(u)\,\widehat{t}_{n-1}(S|u)\;du\nonumber\\
&&\hspace*{1.8cm}+{\beta\over\alpha}\,\int_{r_1}^Sk(u)\,\widehat{t}_{n-1}(S|u)\;du\Biggr\}\qquad(n=1,2,\ldots;\;x<S),
\label{eq:momentsFET}
\end{eqnarray}
where $\widehat{t}_0(S|x)={\widehat P}(S|x)=1$.
\end{theorem}
\begin{proof}
A derivation of (\ref{eq:momentsFET}) follows from the properties of elastic boundaries. We consider again 
separately the cases {\it (i)} and {\it (ii)}. 
\par\noindent
{\it Case (i)}\quad Let $\widehat{t}_n(S_1,S|x)$ $(r_1<S_1<x<S)$ be the first exit time moments  
through the elastic boundary $S$ in the presence of an absorbing boundary  $S_1$. This  
is solution of the differential equation
\begin{equation}
A_1(x)\,{d\psi_n(x)\over dx}+{A_2(x)\over 2}\,{d^2\psi_n(x)\over dx^2}=-n\,\psi_{n-1}(x)
\label{eq:diffequationmoments}
\end{equation}
subject to conditions
\begin{equation}
\lim_{x\downarrow S_1} \psi_n(x)=0,\qquad
 \alpha\,\lim_{x\uparrow S} \psi_n(x)+\beta\,\lim_{x\uparrow S}
\Bigl\{h^{-1}(x)\,{d\psi_n(x)\over dx}\Bigr\}=0.
\label{eq:diffequationmomcond1}
\end{equation}
The general solution of  (\ref{eq:diffequationmoments}) is
\begin{equation}
\psi_n(x)=A+B\int^xh(z)\;dz-n\,\int^xh(z)\;dz\int^zk(u)\,\psi_{n-1}(u)\;du,
\label{eq:gensolutionmom}
\end{equation}
where $A$ and $B$ are arbitrary real constants. By imposing boundary conditions 
(\ref{eq:diffequationmomcond1}), one has
\begin{eqnarray}
&&\widehat{t}_n(S_1,S|x)={n\over\alpha\,{\ds\int_{S_1}^Sh(z)\;dz}+\beta}\;
\Biggl\{\beta\,\int_{S_1}^xh(z)\;dz\int_z^Sk(y)\;\widehat{t}_{n-1}(S_1,S|y)\;dy\nonumber\\
&&\hspace*{2.5cm} +\alpha\;\biggl[\,\int_{S_1}^xh(u)\;du\,\cdot\int_x^Sh(z)\;dz\int_x^zk(y)\;
\widehat{t}_{n-1}(S_1,S|y)\;dy\nonumber\\
&&\hspace*{3.3cm}+\int_x^Sh(u)\;du\,\cdot\int_{S_1}^xh(z)\;dz\int_z^xk(y)\;
\widehat{t}_{n-1}(S_1,S|y)\;dy\biggr]\Biggr\}.
\label{eq:momentsFETdim}
\end{eqnarray}
Since $r_1$ is a natural nonattracting boundary  and $K(r_1,y]<+\infty$, making use of 
(\ref{eq:momentsFETdim}) and by applying l'Hospital's rule,  one has
\begin{eqnarray*}
&&\hspace*{-0.5cm}\widehat{t}_n(S|x):=\lim_{S_1\downarrow r_1}\widehat{t}_n(S_1,S|x)
=n\;\Biggl\{ {\beta\over\alpha}\;
\int_{r_1}^Sk(y)\;\widehat{t}_{n-1}(S|y)\;dy\nonumber\\
&&\hspace*{1.2cm}+\int_x^Sh(z)\;dz\int_x^zk(y)\;\widehat{t}_{n-1}(S_1,S|y)\;dy
+\int_x^Sh(u)\;du\,\cdot\int_{r_1}^xk(y)\;\widehat{t}_{n-1}(S_1,S|y)\;dy\Biggr\},
\end{eqnarray*}
that identifies with  the the right-hand side of (\ref{eq:momentsFET}).
\par\noindent
{\it Case (ii)}\quad Let $\widehat{t}_n(S|x)$ $(r_1<x<S)$ be the first exit time moments  
through the elastic boundary  $S$ in the presence of a reflecting boundary 
or of an entrance boundary $r_1$. This is solution of the differential equation 
(\ref{eq:diffequationmoments}) subject to conditions
\begin{equation}
\lim_{x\downarrow r_1} h^{-1}(x)\,{d\psi_n(x)\over dx}=0,\qquad
 \alpha\,\lim_{x\uparrow S} \psi_n(x)+\beta\,\lim_{x\uparrow S}
\Bigl\{h^{-1}(x)\,{d\psi_n(x)\over dx}\Bigr\}=0.
\label{eq:diffequationmomcond2}
\end{equation}
From (\ref{eq:diffequationmoments}) one obtains again the general solution  
(\ref{eq:gensolutionmom}). By imposing  boundary conditions 
(\ref{eq:diffequationmomcond2}), one finally is led to (\ref{eq:momentsFET}).
\end{proof}
Note that in the absence of refractoriness,   (\ref{eq:momentsFET}) are in agreement 
with (\ref{eq:momentsFPT}). Indeed, if $\beta=0$ one has $\widehat{t}_n(S|x)= t_n(S|x)$. 
The following remark shows that FET moments $\widehat{t}_n(S|x)$ are related to  
the FPT moments $t_n(S|x)$. 
\par
%***************
% REMARK 2.1
%***************
\begin{remark}
Under the assumption {\it (i)} and {\it (ii)},  if $\alpha>0$   one has
\begin{eqnarray}
&&\widehat{t}_n(S|x)=t_n(S|x)+n\,{\beta\over\alpha}\,\sum_{j=0}^{n-1}{n-1\choose j}\,t_{n-1-j}(S|x)\,
\int_{r_1}^Sk(u)\;\widehat{t}_j(S|u)\;du,
\label{eq:relmoments1}\\
&&\widehat{t}_n(S|x)=t_n(S|x)+n\,{\beta\over\alpha}\,\sum_{j=0}^{n-1}{n-1\choose j}\,\widehat{t}_j(S|x)\,
\int_{r_1}^Sk(u)\;t_{n-1-j}(S|u)\;du.
\label{eq:relmoments2}
\end{eqnarray}
\end{remark}
\begin{proof}
Making use of (\ref{eq:momentsFET}),  relations (\ref{eq:relmoments1}) and (\ref{eq:relmoments2}) immediately 
follow by induction.
\end{proof}
\par
Setting $n=1$ in (\ref{eq:relmoments1}) or in (\ref{eq:relmoments2}),  one can see that the mean 
of first exit time is given by
\begin{equation}
\widehat{t}_1(S|x)=t_1(S|x)+{\beta\over\alpha}\,\int_{r_1}^Sk(u)\;du\qquad(x<S).
\label{eq:meanFET}
\end{equation}
Furthermore, setting $n=2$ in (\ref{eq:relmoments1}) and in (\ref{eq:relmoments2}), one can obtain  
two equivalent expressions for the second order moment of first exit time:
\begin{eqnarray}
&&\widehat{t}_2(S|x)=t_2(S|x)+2\;{\beta\over\alpha}\;\widehat{t}_1(S|x)\,\int_{r_1}^Sk(u)\;du
+2\;{\beta\over\alpha}\;\int_{r_1}^Sk(u)\,t_1(S|u)\;du\nonumber\\
&&\hspace*{1.3cm} =t_2(S|x)+2\;{\beta\over\alpha}\;t_1(S|x)\,\int_{r_1}^Sk(u)\;du
+2\;{\beta\over\alpha}\;\int_{r_1}^Sk(u)\,\widehat{t}_1(S|u)\;du.
\label{eq:secondmomentFET}
\end{eqnarray}
Hence, the variance $\widehat{V}(S|x)$ of the first exit time is given by
\begin{eqnarray}
&&\widehat{V}(S|x)=V(S|x)+\biggl({\beta\over\alpha}\,\int_{r_1}^Sk(u)\;du\biggr)^2
+2\;{\beta\over\alpha}\;\int_{r_1}^Sk(u)\,t_1(S|u)\;du,
\label{eq:varianceFET}
\end{eqnarray}
where $V(S|x)$ denotes the FPT variance.
%
%***************
% THEOREM 2.3
%***************
\begin{theorem}
Under the assumption {\it (i)} and {\it (ii)},  if $\alpha>0$ the refractoriness  
period is doomed to end  with certainty and  its moments can be iteratively calculated as
\begin{equation}
E(T_r^n):=\int_0^{\infty}t^n\,g_r(S,t|S)\;dt=n\,{\beta\over\alpha}\,\int_{r_1}^Sk(u)\,
\widehat{t}_{n-1}(S|u)\;du\qquad(n=1,2,\ldots).
\label{eq:momentsrefractoriness}
\end{equation}
\end{theorem}
\begin{proof}
Integrating  both sides of (\ref{eq:integralequat}) in $(0,+\infty)$ one has  
$$
\int_0^{+\infty}g_r(S,t|S)\;dt=1,
$$
implying that  the refractoriness  period is doomed to end  with certainty. Furthermore, 
from (\ref{eq:integralequat})  we also have: 
\begin{eqnarray*}
&&\widehat{t}_n(S|x):=\int_0^{+\infty}t^n\,g_e(S,t|x)\;dt
=\int_0^{+\infty}dt\;t^n\;\int_0^tg(S,\tau|x)\,g_r(S,t|S,\tau)\;d\tau\nonumber\\
&&\hspace*{1.3cm}=\int_0^{+\infty}d\tau\;g(S,\tau|x)\;\int_{\tau}^{+\infty}t^n\,
g_r(S,t|S,\tau)\;dt\nonumber\\
&&\hspace*{1.3cm}=\int_0^{+\infty}d\tau\;g(S,\tau|x)\;
\int_0^{+\infty}(\tau+\vartheta)^ng_r(S,\vartheta|S)\;d\vartheta\nonumber\\
&&\hspace*{1.3cm}=\sum_{j=0}^n\;{n\choose j}\,t_{n-j}(S|x)\,E\bigl(T_r^j\bigr).\qquad (n=1,2,\ldots).
\label{eq:relationmoments1}
\end{eqnarray*}
Hence, 
\begin{equation}
E\bigl(T_r^n\bigr)=\widehat{t}_n(S|x)-\sum_{j=0}^{n-1}\;{n\choose j}\,t_{n-j}(S|x)\,E\bigl(T_r^j
\bigr)\qquad(n=1,2,\ldots).
\label{eq:relationmoments2}
\end{equation}
We now proceed by induction. Setting $n=1$ in (\ref{eq:relationmoments2}) one sees  
that $E\bigl(T_r\bigr)=\widehat{t}_1(S|x)-t_1(S|x)$. Hence, 
on account of (\ref{eq:meanFET}), (\ref{eq:momentsrefractoriness}) holds for $n=1$.   
Furthermore, assuming that  (\ref{eq:relationmoments2}) hold for $j=1,2,\ldots,n$, the right-hand side of 
(\ref{eq:relationmoments2}) for $n+1$ becomes:
\begin{eqnarray}
&&\hspace*{-0.4cm}\widehat{t}_{n+1}(S|x)-\sum_{j=0}^n\;{n+1\choose j}\,
t_{n+1-j}(S|x)\,E\bigl(T_r^j\bigr)\nonumber\\
&&\hspace*{1.2cm} =\widehat{t}_{n+1}(S|x)-t_{n+1}(S|x)-(n+1)\,{\beta\over\alpha}\,
\sum_{j=0}^{n-1}{n\choose j}\,t_{n-j}(S|x)\int_{r_1}^Sk(u)\;\widehat{t}_j(S|u)\;du\nonumber\\
&&\hspace*{1.2cm} =(n+1)\,{\beta\over\alpha}\,\int_{r_1}^Sk(u)\;
\widehat{t}_n(S|u)\;du,
\label{eq:relationmoments3}
\end{eqnarray}
where the last equality follows from (\ref{eq:relmoments1}). From (\ref{eq:relationmoments2}) we 
note that the left-hand side of (\ref{eq:relationmoments3}) is equal to $E\bigl(T_r^{n+1}\bigr)$. 
Hence,  if (\ref{eq:momentsrefractoriness}) holds for an arbitrarily fixed $n$, it also holds 
for $n+1$, which completes the proof. 
\end{proof}
Comparing  (\ref{eq:momentsFET}) and (\ref{eq:momentsrefractoriness}) we note that
$$
E(T_r^n)\equiv \lim_{x\uparrow S}\;\widehat{t}_n(S|x).
$$
In particular, from (\ref{eq:momentsrefractoriness}) the first two moments and the variance of 
the refractoriness  period are seen to be:
\begin{eqnarray}
&& E(T_r)={\beta\over\alpha}\,\int_{r_1}^Sk(u)\;du\nonumber\\
&& E(T_r^2)=2\;{\beta\over\alpha}\;\int_{r_1}^Sk(u)\;t_1(S|u)\;du
+2\,\biggl({\beta\over\alpha}\,\int_{r_1}^Sk(u)\;du\biggr)^2
\label{eq:firstmoments}\\
&& V(T_r)=2\;{\beta\over\alpha}\;\int_{r_1}^Sk(u)\;t_1(S|u)\;du
+\biggl({\beta\over\alpha}\,\int_{r_1}^Sk(u)\;du\biggr)^2.\nonumber
\end{eqnarray}
Comparing the first and last of (\ref{eq:firstmoments}) with (\ref{eq:meanFET}) and 
(\ref{eq:varianceFET}), we have
\begin{equation}
\widehat{t}_1(S|x)=t_1(S|x)+E(T_r),\qquad \widehat{V}(S|x)=V(S|x)+V(T_r),
\label{eq:relfirstmoments}
\end{equation}
i.e. the mean (variance) of first exit time through $S$ starting from $x$ 
is the sum of the mean (variance) of first passage time through $S$ 
starting from $x$ and of the mean (variance) of the refractoriness period.
%------------------------------------------------------
\section{Analysis of three neuronal models}
%------------------------------------------------------
\hfill\\
In order to embody some physiological features of real neurons, several 
alternative models have been proposed in the literature (cf, for instance, 
\cite{Ricciardi99}, \cite{Ricciardi02} and references therein). In this Section 
we shall  investigate the behavior of the refractoriness  period 
for the Wiener,  Ornstein-Uhlenbeck (OU) and Feller neuronal models. 
We assume that all three  neuronal models are restricted to the same diffusion interval 
$I=[\nu,+\infty)$, having set $r_1=\nu$.
%
% WIENER MODEL
%======================= TABLE 1 Wiener model ======================================
\begin{table}[h]
 \begin{center}
{\footnotesize
 \begin{tabular}{rlllll}
 \hline
 & $t_1(S|\rho)$ & $E(T_r)$ & $E(T_r)$  & $E(T_r)$  & $E(T_r)$\\
$\sigma^2$ & & $p_{_{R}}=0.1$ & $p_{_{R}}=0.5$  & $p_{_{R}}=0.9$  & $p_{_{R}}=0.99$\\
\hline
10.  & 3.073451 E+2  & 6.294544 E+2    & 5.665090 E+3 & 5.098581 E+4   & 5.608439 E+5\\
20.  & 7.331871 E+1  & 9.425701 E+0    & 8.483131 E+1 & 7.634818 E+2   & 8.398300 E+3\\
30.  & 3.936016 E+1  & 2.021650 E+0    & 1.819485 E+1 & 1.637537 E+2   & 1.801290 E+3\\
40.  & 2.663797 E+1  & 8.663807 E$-1$  & 7.797426 E+0 & 7.017684 E+1   & 7.719452 E+2\\
50.  & 2.007160 E+1  & 4.966112 E$-1$  & 4.469501 E+0 & 4.022551 E+1   & 4.424806 E+2\\
100. & 8.937578 E+0  & 1.281821 E$-1$  & 1.153639 E+0 & 1.038275 E+1   & 1.142103 E+2\\
200. & 4.225259 E+0  & 4.617762 E$-2$  & 4.155986 E$-1$ & 3.740387 E+0 & 4.114426 E+1\\
300. & 2.765483 E+0  & 2.760995 E$-2$  & 2.484895 E$-1$ & 2.236406 E+0 & 2.460046 E+1\\
400. & 2.055224 E+0  & 1.961207 E$-2$  & 1.765086 E$-1$ & 1.588577 E+0 & 1.747435 E+1\\
500. & 1.635207 E+0  & 1.518666 E$-2$  & 1.366799 E$-1$ & 1.230119 E+0 & 1.353131 E+1\\
\hline
\end{tabular}
}
\end{center}
\caption{{\small Wiener model with $\mu=-0.5$ and 
$\sigma^2=10\cdot i,  100\cdot i\;(i=1,2,\ldots,5)$, 
restricted to  $I=[\nu,+\infty)$ with $\nu=-80$. In the second column we have listed the FPT 
mean $t_1(S|\rho)$ with $S=-50$ and $\varrho=-70$. Instead, in columns three, four, five and six 
we have respectively listed the mean of refractoriness period for  $p_{_{R}}=0.1,\,0.5,\,0.9,\,0.99$.}}
\label{table1}
\end{table} 
%========================================================================
%
%======================= TABLE 2 Wiener model ======================================
\begin{table}[ht]
 \begin{center}
{\footnotesize
 \begin{tabular}{rlllll}
 \hline
 & $V(S|\rho)$ & $V(T_r)$ & $V(T_r)$  & $V(T_r)$  & $V(T_r)$\\
$\sigma^2$ &  & $p_{_{R}}=0.1$ & $p_{_{R}}=0.5$  & $p_{_{R}}=0.9$  & $p_{_{R}}=0.99$\\
\hline
10.  & 9.254218 E+4 & 7.681238 E+5   & 3.544044 E+7    & 2.629677 E+9  & 3.148772 E+11\\
20.  & 4.970295 E+3 & 1.310444 E+3   & 1.819075 E+4    & 6.818541 E+5  & 7.161989 E+7\\
30.  & 1.390880 E+3 & 1.385660 E+2   & 1.541363 E+3    & 3.770806 E+4  & 3.364468 E+6\\
40.  & 6.265060 E+2 & 3.874901 E+1   & 4.027854 E+2    & 8.002658 E+3  & 6.297560 E+5\\
50.  & 3.519348 E+2 & 1.638408 E+1   & 1.652135 E+2    & 2.925226 E+3  & 2.101676 E+5\\
100. & 6.821593 E+1 & 1.804893 E+0   & 1.742704 E+1    & 2.526670 E+2  & 1.463751 E+4\\
200. & 1.506377 E+1 & 3.008332 E$-1$ & 2.861029 E+0    & 3.818526 E+1  & 1.958992 E+3\\
300. & 6.426684 E+0 & 1.168741 E$-1$ & 1.106754 E+0    & 1.440657 E+1  & 7.086384 E+2\\
400. & 3.542131 E+0 & 6.147195 E$-2$ & 5.809411 E$-1$  & 7.471651 E+0  & 3.597818 E+2\\
500. & 2.239493 E+0 & 3.778974 E$-2$ & 3.567134 E$-1$  & 4.555481 E+0  & 2.165615 E+2\\
\hline
\end{tabular}
}
\end{center}
\caption{{\small For the Wiener model and for the same choices of parameters 
of Table \ref{table1},  in the second column we have listed the FPT variance 
$V(S|\varrho)$ with $S=-50$ and $\varrho=-70$, whereas in columns three, 
four, five and six we have respectively listed 
the variance of refractoriness period for  $p_{_{R}}=0.1,\,0.5,\,0.9,\,0.99$.}}
\label{table2}
\end{table} 
%========================================================================
%
%------------------------------------------------------
\subsection{Wiener model}
%------------------------------------------------------
The Wiener neuronal model is defined as the diffusion process 
$X(t)$ characterized by the following drift and infinitesimal variance: 
\begin{equation}
A_1(x)=\mu\qquad\qquad A_2(x)=\sigma^2,\qquad (\mu\in\Real,\sigma>0),
\label{eq:Wienermodel}
\end{equation}
restricted to $I=[\nu,+\infty)$, where on the regular boundary $x=\nu$ 
a reflecting condition is imposed.  For such 
process  scale and speed functions are
$$
h(x)=\exp\Bigl\{-{2\,\mu\, x\over\sigma^2}\Bigr\},\qquad 
k(x)={2\over\sigma^2}\,\exp\Bigl\{{2\,\mu\, x\over\sigma^2}\Bigr\}.
$$
Furthermore, the mean of first passage time is
\begin{equation}
t_1(S|x)=\left\{ \begin{array}{ll}
{\ds{(S-x)\,(S+x-2\,\nu)\over\sigma^2}}, &\mu=0 \\
\\
{\ds{S-x\over\mu}+{\sigma^2\over 2\,\mu^2}
\biggl[\exp\Bigl\{-\,{2\,\mu\,(S-\nu)\over\sigma^2}\Bigr\}-
\exp\Bigl\{-\,{2\,\mu\,(x-\nu)\over\sigma^2}\Bigr\}\biggr]}\,, &\mu\neq 0.
\end{array}
\right.
\label{eq:FPTmeanWiener}
\end{equation}
For the Wiener model (\ref{eq:Wienermodel}) with $\mu=-0.5$, 
$\sigma^2=10\cdot i, 100\cdot i \;(i=1,2,\ldots,5)$, restricted 
to  $I=[\nu,+\infty)$ with $r_1\equiv\nu=-80$, in the second 
column of Table \ref{table1}  and of Table \ref{table2} we have respectively 
listed the mean $t_1(S|\varrho)$  and variance  $V(S|\varrho)$,  
numerically obtained via (\ref{eq:momentsFPT}) with $S=-50$ and  $\varrho=-70$. 
Note that the FPT mean and variance decrease with $\sigma^2$. Being 
$\beta/\alpha=p_{_{R}}/(1-p_{_{R}})$, 
in Table \ref{table1} and in Table \ref{table2} we have 
respectively listed the values of mean and variance of refractoriness period, 
numerically obtained via  (\ref{eq:firstmoments}) for $p_{_{R}}=0.1,\,0.5,\,0.9,\,0.99$. 
We observe that $E(T_r)$ and  $V(T_r)$ increase with $p_{_{R}}$ for any fixed $\sigma^2$.
%
% OU MODEL
%
%======================= TABLE 3 OU model ======================================
\begin{table}[ht]
 \begin{center}
{\footnotesize
 \begin{tabular}{rlllll}
 \hline
& $t_1(S|\rho)$ & $E(T_r)$ & $E(T_r)$  & $E(T_r)$  & $E(T_r)$\\
$\sigma^2$ & & $p_{_{R}}=0.1$ & $p_{_{R}}=0.5$  & $p_{_{R}}=0.9$  & $p_{_{R}}=0.99$\\
\hline
10.  & 9.862135 E+3 & 9.901436 E+41   & 8.911293 E+42 & 8.020163 E+43 & 8.822180 E+44\\
20.  & 2.600359 E+2 & 3.452097 E+20   & 3.106887 E+21 & 2.796199 E+22 & 3.075818 E+23\\
30.  & 7.956655 E+1 & 2.140293 E+13   & 1.926264 E+14 & 1.733637 E+15 & 1.907001 E+16\\
40.  & 4.273886 E+1 & 4.978530 E+9    & 4.480677 E+10 & 4.032609 E+11 & 4.435870 E+12\\
50.  & 2.853092 E+1 & 3.149993 E+7    & 2.834994 E+8  & 2.551494 E+9  & 2.806644 E+10\\
100. & 1.038152 E+1 & 1.006196 E+3    & 9.055763 E+3  & 8.150187 E+4  & 8.965206 E+5\\
200. & 4.525217 E+0 & 4.073683 E+0    & 3.666314 E+1  & 3.299683 E+2  & 3.629651 E+3\\
300. & 2.890905 E+0 & 5.465816 E$-1$  & 4.919234 E+0  & 4.427311 E+1  & 4.870042 E+2\\
400. & 2.123662 E+0 & 1.839895 E$-1$  & 1.655905 E+0  & 1.490315 E+1  & 1.639346 E+2\\
500. & 1.678216 E+0 & 9.103197 E$-2$  & 8.192877 E$-1$  & 7.373589 E+0  & 8.110948 E+1\\
\hline
\end{tabular}
}
\end{center}
\caption{{\small OU model with $\vartheta=5$, $\varrho=-70$ and  
$\sigma^2=10\cdot i, 100\cdot i \;(i=1,2,\ldots,5)$, 
restricted to  $I=[\nu,+\infty)$ with $\nu=-80$. 
In the second column  we have listed the FPT mean $t_1(S|\rho)$ with $S=-50$, whereas in columns three, 
four, five and six we have respectively listed the mean of refractoriness period 
for  $p_{_{R}}=0.1,\,0.5,\,0.9,\,0.99$.}}
\label{table3}
\end{table} 
%========================================================================
%
%======================= TABLE 4 OU model ======================================
\begin{table}[ht]
 \begin{center}
{\footnotesize
 \begin{tabular}{rlllll}
 \hline
& $V(S|\rho)$ & $V(T_r)$ & $V(T_r)$  & $V(T_r)$  & $V(T_r)$\\
$\sigma^2$ & & $p_{_{R}}=0.1$ & $p_{_{R}}=0.5$  & $p_{_{R}}=0.9$  & $p_{_{R}}=0.99$\\
\hline
10.  & 9.713857 E+7   & 9.803844 E+83   & 7.941114 E+85  & 6.432302 E+87  & 7.783086 E+89\\
20.  & 6.554937 E+4   & 1.191697 E+41   & 9.652749 E+42  & 7.818727 E+44  & 9.460659 E+46\\
30.  & 5.898427 E+3   & 4.580854 E+26   & 3.710492 E+28  & 3.005498 E+30  & 3.636653 E+32\\
40.  & 1.654790 E+3   & 2.478576 E+19   & 2.007647 E+21  & 1.626194 E+23  & 1.967694 E+25\\
50.  & 7.239524 E+2   & 9.922473 E+14   & 8.037191 E+16  & 6.510124 E+18  & 7.877250 E+20\\
100. & 9.240940 E+1   & 1.029849 E+6    & 8.216362 E+7   & 6.643966 E+9   & 8.037646 E+11\\
200. & 1.728177 E+1   & 4.589508 E+1    & 1.607888 E+3   & 1.112524 E+5   & 1.320047 E+7\\
300. & 7.020488 E+0   & 2.765114 E+0    & 4.639613 E+1   & 2.159884 E+3   & 2.393706 E+5\\
400. & 3.780330 E+0   & 6.379705 E$-1$  & 8.179088 E+0   & 2.710374 E+2   & 2.741284 E+4\\
500. & 2.357829 E+0   & 2.431146 E$-1$  & 2.784683 E+0   & 7.339088 E+1   & 6.787980 E+3\\
\hline
\end{tabular}
}
\end{center}
\caption{{\small For the OU model and for the same choices of parameters 
of Table \ref{table4},  in the second column we have listed the FPT variance 
$V(S|\varrho)$ with $S=-50$ and $\varrho=-70$, 
whereas in columns three, four, five and six we have respectively listed 
the variance of refractoriness period for  $p_{_{R}}=0.1,\,0.5,\,0.9,\,0.99$.}}
\label{table4}
\end{table} 
%========================================================================
%
%------------------------------------------------------
\subsection{OU model}
%------------------------------------------------------
The OU neuronal model is defined as the diffusion process 
$X(t)$ characterized by the following drift and infinitesimal variance: 
\begin{equation}
A_1(x)=-{1\over\vartheta}\,(x-\varrho),\qquad A_2=\sigma^2
\qquad(\varrho\in {\bf R},\sigma>0,\,\vartheta>0), 
\label{eq:OUmodel}
\end{equation}
restricted to $I=[\nu,+\infty)$, where on the regular boundary $x=\nu$ 
a reflecting condition is imposed.  For such process  the scale and speed functions are
$$
h(x)=\exp\Bigl\{ {x^2\over\vartheta\,\sigma^2}
-{2\,\varrho\,x\over\vartheta\,\sigma^2}\Bigr\},\qquad 
k(x)={2\over\sigma^2}\,\exp\Bigl\{ -{x^2\over\vartheta\,\sigma^2}
+{2\,\varrho\,x\over\vartheta\sigma^2}\Bigr\}.
$$
Furthermore, the mean of first passage time is:
\begin{eqnarray}
&& t_1(S|x)= \vartheta\,\sum_{k=0}^{+\infty}{2^k\over(k+1)\,(2k+1)!!}\,
\biggl[\biggr({S-\varrho\over\sigma\sqrt{\vartheta}}\biggr)^{2k+2}
-\biggr({x-\varrho\over\sigma\sqrt{\vartheta}}\biggr)^{2k+2}\biggr]\nonumber\\
&&\hspace*{1.5cm} -2\,\vartheta\,
\exp\Bigl\{-\,{(\nu-\varrho)^2\over\sigma^2\vartheta}\Bigr\}\,
\sum_{k=0}^{+\infty}{2^k\over(2k+1)!!}\,
\biggr({\nu-\varrho\over\sigma\sqrt{\vartheta}}\biggr)^{2k+1}\nonumber\\
&&\hspace*{1.5cm}\times\sum_{k=0}^{+\infty}{1\over(2k+1)\,k!}\,
\biggl[\biggr({S-\varrho\over\sigma\sqrt{\vartheta}}\biggr)^{2k+1}
-\biggr({x-\varrho\over\sigma\sqrt{\vartheta}}\biggr)^{2k+1}\biggr].
\label{eq:FPTmeanOU}
\end{eqnarray}
For the OU model (\ref{eq:OUmodel}) with $\vartheta=5$, $\varrho=-70$, 
$\sigma^2=10\cdot i, 100\cdot i \;(i=1,2,\ldots,5)$, restricted to $I=[\nu,+\infty)$ 
with $\nu=-80$, in the second column of Table \ref{table3}  and of Table \ref{table4} 
we have respectively listed the mean $t_1(S|\varrho)$  and variance  $V(S|\varrho)$,  
numerically obtained via (\ref{eq:momentsFPT}) with $S=-50$. 
 Furthermore,  
in Table \ref{table3} and in Table \ref{table4} we have also 
listed the values of mean and variance of refractoriness period for $p_{_{R}}=0.1,\,0.5,\,0.9,\,0.99$. 
Similarly to the case of the Wiener model, for the OU model the FPT mean and variance decrease 
with $\sigma^2$; furthermore, $E(T_r)$ and  $V(T_r)$ increase with $p_{_{R}}$ for any fixed $\sigma^2$.
%
% NEW FELLER MODEL
%
%======================= TABLE 5 Feller model ======================================
\begin{table}[h]
 \begin{center}
{\footnotesize
 \begin{tabular}{rlllll}
 \hline
& $t_1(S|\rho)$ & $E(T_r)$ & $E(T_r)$  & $E(T_r)$  & $E(T_r)$\\
$\xi$ & & $p_{_{R}}=0.1$ & $p_{_{R}}=0.5$  & $p_{_{R}}=0.9$  & $p_{_{R}}=0.99$\\
\hline
0.5 & 3.768002 E+2 & 4.103229 E+15  & 3.692906 E+16  & 3.323615 E+17 & 3.655977 E+18\\
1.0 & 8.029989 E+1 & 2.425535 E+7   & 2.182981 E+8   & 1.964683 E+9  & 2.161152 E+10\\
1.5 & 4.661249 E+1 & 4.020549 E+4   & 3.618494 E+5   & 3.256645á E+6 & 3.582309 E+7\\
2.0 & 3.470051 E+1 & 1.573636 E+3   & 1.416272 E+4   & 1.274645 E+5  & 1.402110  E+6\\
2.5 & 2.866867 E+1 & 2.204449 E+2   & 1.984004 E+3   & 1.785603 E+4  & 1.964164 E+5\\
3.0 & 2.502681 E+1 & 5.871921 E+1   & 5.284729 E+2   & 4.756256 E+3  & 5.231882 E+4\\
3.5 & 2.258692 E+1 & 2.264293 E+1   & 2.037863 E+2   & 1.834077 E+3  & 2.017485 E+4\\
4.0 & 2.083633 E+1 & 1.102071 E+1   & 9.918635 E+1   & 8.926772 E+2  & 9.819449 E+3\\
4.5 & 1.951789 E+1 & 6.271095 E+0   & 5.643986 E+1   & 5.079587 E+2  & 5.587546 E+3\\
5.0 & 1.848842 E+1 & 3.983514 E+0   & 3.585162 E+1   & 3.226646 E+2  & 3.549311 E+3\\
 \hline
\end{tabular}
}
\end{center}
\caption{{\small Feller model with $\vartheta=5$, $\varrho=-70$,  $\nu=-80$ and  
$\xi= 0.5\cdot i\;(i=1,2,\ldots, 10)$. In second column  we have listed the 
FPT mean $t_1(S|\rho)$ with $S=-50$, whereas in columns three, 
four, five and six we have respectively listed the mean of refractoriness period 
for  $p_{_{R}}=0.1,\,0.5,\,0.9,\,0.99$.}}
\label{table5}
\end{table} 
%========================================================================
%
%======================= TABLE 6 Feller model ======================================
\begin{table}[h]
 \begin{center}
{\footnotesize
 \begin{tabular}{rlllll}
 \hline
& $V(S|\rho)$ & $V(T_r)$ & $V(T_r)$  & $V(T_r)$  & $V(T_r)$\\
$\xi$ & & $p_{_{R}}=0.1$ & $p_{_{R}}=0.5$  & $p_{_{R}}=0.9$  & $p_{_{R}}=0.99$\\
\hline
0.5 & 1.395404 E+5 & 1.683649 E+31  & 1.363755 E+33 & 1.104642 E+35 & 1.336617 E+37\\
1.0 & 6.372482 E+3 & 5.883257 E+14  & 4.765411 E+16 & 3.859980 E+18 & 4.670576 E+20\\
1.5 & 2.241795 E+3 & 1.620153 E+9   & 1.309681 E+11 & 1.060603 E+13 & 1.283297 E+15\\
2.0 & 1.304116 E+3 & 2.585121 E+6   & 2.015619 E+8  & 1.625601 E+10 & 1.966008 E+12\\
2.5 & 9.313963 E+2 & 6.143569 E+4   & 4.051826 E+6  & 3.198777 E+8  & 3.85908 E+10\\
3.0 & 7.390905 E+2 & 6.492913 E+3   & 3.066862 E+5  & 2.286843 E+7  & 2.739949 E+9\\
3.5 & 6.238662 E+2 & 1.592475 E+3   & 5.124507 E+4  & 3.451142 E+6  & 4.079672 E+8\\
4.0 & 5.478171 E+2 & 6.147816 E+2   & 1.427634 E+4  & 8.366949 E+5  & 9.684436 E+7\\
4.5 & 4.940875 E+2 & 3.064866 E+2   & 5.588538 E+3  & 2.795397 E+5  & 3.144373 E+7\\
5.0 & 4.541290 E+2 & 1.789263 E+2   & 2.751622 E+3  & 1.172086 E+5  & 1.272925 E+7\\
\hline
\end{tabular}
}
\end{center}
\caption{{\small For the Feller model and for the same choices of parameters 
of Table \ref{table5},  in the second column we have listed the FPT variance $V(S|\varrho)$ 
with $S=-50$ and $\varrho=-70$, whereas in columns three, four, five and six we have 
respectively listed the variance of refractoriness period for  $p_{_{R}}=0.1,\,0.5,\,0.9,\,0.99$.}}
\label{table6}
\end{table} 
%========================================================================
%------------------------------------------------------
\subsection{Feller model}
%------------------------------------------------------
The Feller neuronal model is defined as the diffusion process 
$X(t)$ characterized by the following drift and infinitesimal variance: 
\begin{equation}
A_1(x)=-{1\over\vartheta}\,(x-\varrho),\qquad A_2(x)=2\,\xi\,(x-\nu)\qquad
(\varrho,\nu\in{\bf R},\varrho>\nu,\vartheta>0,\xi>0).
\label{eq:Fellermodel}
\end{equation}
defined in $I=[\nu,+\infty)$, where  $x=\nu$ is regular 
if $\varrho-\nu<\xi\,\vartheta$ and entrance if $\varrho-\nu\geq\xi\,\vartheta$,  
whereas the  boundary $+\infty$ is natural. For such process  the scale 
and speed functions are:
$$
h(x)=\exp\Bigl\{ {x\over\vartheta\,\xi}\Bigr\}\,
\bigl(x-\nu\bigr)^{-(\varrho-\nu)/(\vartheta\,\xi)},\qquad 
k(x)={1\over\xi}\,\exp\Bigl\{-{x\over\vartheta\,\xi}\Bigr\}\,
\bigl(x-\nu\bigr)^{(\varrho-\nu)/(\vartheta\,\xi)-1}.
$$
The mean of the firing time can be calculated;  for $x<S$ one obtains
\begin{equation}
t_1(S|x)={\vartheta\over \varrho-\nu}\,\biggl[S-x+\sum_{k=1}^\infty\biggl({1\over\vartheta}\biggr)^k
{(S-\nu)^{k+1}-(x-\nu)^{k+1}\over k+1}\;
\biggl\{\prod_{i=1}^k\biggl({\varrho-\nu\over\vartheta}+\xi\,i\biggr)\biggr\}^{-1}\biggr].
\label{eq:tm_Feller}
\end{equation}
For the Feller model (\ref{eq:Wienermodel}) with $\vartheta=5$, $\varrho=-70$, $\nu=-80$, 
$\xi=0.5\cdot i \;(i=1,2,\ldots,10)$, in the second 
column of Table \ref{table5}  and of Table \ref{table6} we have respectively 
listed the mean $t_1(S|\varrho)$  and variance  $V(S|\varrho)$,  
numerically obtained via (\ref{eq:momentsFPT}) with $S=-50$. 
Note that the FPT mean and variance decrease with $\xi$. Furthermore, 
in Table \ref{table5} and in Table \ref{table6} we have  listed 
the values of mean and variance of refractoriness period for $p_{_{R}}=0.1,\,0.5,\,0.9,\,0.99$. 
We note that $E(T_r)$ and  $V(T_r)$ increase with $p_{_{R}}$ for any fixed $\xi$.
\par
We conclude by pointing out that the purpose of the present note was to establish the 
quantitative foundations to a viable way to include refractoriness in neuronal 
diffusion models. Implementation of our approach to data analysis will be the object 
of future endeavors.

\ack Work performed within a joint cooperation agreement between
Japan Science and Technology Corporation (JST) and Universit\`a di Napoli
Federico II, under partial support by INdAM (G.N.C.S).
%--------------------------------------------------------------------
%
%REFERENCES
%
%

%
\end{document}